\documentclass[10pt]{amsart}

\usepackage{amssymb,amsmath}
\usepackage[dvips]{graphicx,color}

\newtheorem{theorem}{Theorem}
\newtheorem{lemma}{Lemma}
\newtheorem{example}{Example}
\newtheorem{problem}{Problem}
\newcommand{\bt}{\begin{theorem}}
\newcommand{\et}{\end{theorem}}
\newcommand{\bl}{\begin{lemma}}
\newcommand{\el}{\end{lemma}}
\newcommand{\bex}{\begin{example}}
\newcommand{\eex}{\end{example}}
\newcommand{\bp}{\begin{problem}}
\newcommand{\ep}{\end{problem}}
\newcommand{\beal}{\begin{align*}}
\newcommand{\enal}{\end{align*}}
\newcommand{\beq}{\begin{equation}}
\newcommand{\eeq}{\end{equation}}
\newcommand{\benum}{\begin{enumerate}}
\newcommand{\eenum}{\end{enumerate}}
\newcommand{\ba}{\begin{array}}
\newcommand{\ea}{\end{array}}
\newcommand{\Z}{\ensuremath{\mathbb Z}}
\newcommand{\N}{\ensuremath{\mathbb N}}

\newcommand{\card}{\text{card}}

\begin{document}

\title{A new upper bound for finite additive bases}
\author{C.~S{\. i}nan G{\" u}nt{\" u}rk}
\address{Courant Institute\\New York University\\New York, New York 10012}
\email{gunturk@courant.nyu.edu}
\author{Melvyn B.~Nathanson}
\thanks{The work of S.G. was supported in part by NSF grant DMS 0219072.
The work of M.B.N. was supported in part by grants from the NSA Mathematical
Sciences Program and the PSC-CUNY Research Award Program.}
\address{Department of Mathematics\\Lehman College (CUNY)\\Bronx, New York
10468} \email{melvyn.nathanson@lehman.cuny.edu}

\begin{abstract}
Let $n(2,k)$ denote the largest integer $n$ for which there exists a set
$A$ of $k$ nonnegative integers such that the sumset $2A$ contains
$\{0,1,2,\ldots,n-1\}.$  A classical problem in additive number theory
is to find an upper bound for $n(2,k).$  In this paper it is proved that
$\limsup_{k\rightarrow\infty} n(2,k)/k^2 \leq 0.4789.$
\end{abstract}

\maketitle

\section{An extremal problem for finite bases}
Let $\N_0$ and \Z\ denote the nonnegative integers and integers,
respectively, and let $|A|$ denote the cardinality of the set $A$.

Let $A$ be a set of integers, and consider the {\em sumset}
\[
2A = \{ a+a' : a,a' \in A\}.
\]
Let $S$ be a set of integers.
The set $A$ is a {\em basis of order 2 for $S$} if $S \subseteq 2A$.
The set $A$ is called a {\em basis of order 2 for $n$} if the sumset
$2A$ contains the first $n$ nonnegative integers, that is, if $A$ is a
basis of order 2 for the interval of integers $[0,n-1]:=\{0,1,2,\ldots,n-1\}.$
We define $n(2,A)$ as the largest integer $n$ such that $A$ is a basis of
order 2 for $n$, that is,
\[
n(2,A) = \max\{n: [0,n-1] \subseteq 2A\}.
\]
Rohrbach~\cite{rohr37a} introduced the extremal problem of determining the
largest integer $n$ for which there exists a set $A$ consisting of at most
$k$ nonnegative integers such that $A$ is a basis of order 2 for $n$.  Let
\[
n(2,k) = \max\{n(2,A) : A \subseteq \N_0 \text{ and } |A| = k\}.
\]
Rohrbach's problem is to compute or estimate the extremal function $n(2,k).$
The set $A$ is called an {\em extremal $k$-basis of order 2} if $|A| \leq k$
and $n(2,A) = n(2,k).$

For example, $n(2,1) = 1$ and $n(2,2) = 3$.   The unique extremal 1-basis of
order 2 is $\{0\}$, and the unique extremal 2-basis of order 2 is $\{0,1\}$.
For $k = 3$ we have $n(2,3) = 5,$ and the extremal 3-bases of order 2 are
$\{0,1,2\}$ and $\{0,1,3\}$.
If $k \geq 2$ and $A$ is an extremal $k$-basis of order 2, then $0,1 \in A.$
If $A$ is a finite set of $k$ nonnegative integers and $n(2,A) = n,$ then
$n \not\in A.$  If $a \in A$ and $a > n,$ then the set
$A' = (A\setminus\{a\}) \cup \{n\}$ has cardinality $k$, and
$n(2,A') \geq n+1 > n(2,A).$  Therefore, if $A$ is an extremal $k$-basis
of order 2 and $n(2,k) = n,$ then
\[
\{0,1\} \subseteq A \subseteq \{0,1,2,\ldots,n-1\} \subseteq 2A.
\]
If $A$ is an extremal $k$-basis for $n$, then $|A| = k$ and
$A \subseteq \{0,1,2,\ldots, n-1\}.$

Rohrbach determined order of magnitude of $n(2,k)$.  He observed that if
$A$ is a set of cardinality $k$, then there are exactly  $\binom{k+1}{2}$
ordered pairs of the form $(a,a')$ with $a,a' \in A$ and $a \leq a'$.
This gives the upper bound
\[
n(2,k) \leq \binom{k+1}{2} = \frac{k^2}{2} + O(k).
\]
To derive a lower bound, he set $r = [k/2]$ and constructed the set
\[
A = \{ 0,1,2,\ldots, r-1,r,2r,3r,\ldots,(r-1)r\}.
\]
We have
\[
|A| = 2r-1 \leq k
\]
and $\{0,1,\ldots, r^2\} \subseteq 2A.$
Then
\[
n(2,A) \geq r^2 + 1 \geq \frac{(k-1)^2}{4} + 1 = \frac{k^2}{4} + O(k)
\]
and so
\[
n(2,k) \geq \frac{k^2}{4} + O(k).
\]
Thus,
\[
\liminf_{n\rightarrow\infty}\frac{n(2,k)}{k^2} \geq \frac{1}{4}= 0.25
\]
and
\[
\limsup_{n\rightarrow\infty}\frac{n(2,k)}{k^2} \leq \frac{1}{2} = 0.5.
\]
It is a open problem to compute these upper and lower limits.
Mrose~\cite{mros79,hofm01} proved that
\[
\liminf_{n\rightarrow\infty}\frac{n(2,k)}{k^2} \geq \frac{2}{7} = 0.2857\ldots,
\]
and this is still the best lower bound.
Rohrbach used a combinatorial argument to get the nontrivial upper bound
\[
\limsup_{n\rightarrow\infty}\frac{n(2,k)}{k^2} \leq 0.4992.
\]
Moser~\cite{mose60} introduced Fourier series to obtain
\[
\limsup_{n\rightarrow\infty}\frac{n(2,k)}{k^2} \leq 0.4903,
\]
and subsequent improvements by Moser, Pounder, and
Riddell~\cite{mose-poun-ridd69}
produced
\[
\limsup_{n\rightarrow\infty}\frac{n(2,k)}{k^2} \leq 0.4847.
\]
Combining Moser's analytic method and Rohrbach's combinatorial technique,
Klotz~\cite{klot69b} proved that
\[
\limsup_{n\rightarrow\infty}\frac{n(2,k)}{k^2} \leq 0.4802.
\]
In this paper, we use Fourier series for functions of two variables to obtain
\[
\limsup_{n\rightarrow\infty}\frac{n(2,k)}{k^2} \leq 0.4789.
\]

We note that  Rohrbach used a slightly different function $n(2,k)$:
He defined $n(2,k)$ as the largest integer $n$ for which there
exists a set $A$ consisting of $k+1$ nonnegative integers such that
the sumset $2A$ contains the first $n+1$ nonnegative integers.  Of
course, Rohrbach's function and our function have the same
asymptotics.

\section{Moser's application of Fourier series}
In this section we describe Moser's use of harmonic analysis to obtain an
upper bound for $n(2,k)$. Let $A$ be an extremal $k$-basis of order 2.
Let $r_{2,A}(j)$ denote the number of representations of $j$ as a sum of two
elements of $A$, that is,
\[
r_{2,A}(j) = \card\left( \left\{ (a_1,a_2) \in A \times A : a_1 + a_2 = j
\text{ and } a_1 \leq a_2 \right\}  \right).
\]
We introduce the generating function
\[
f_A(q) = \sum_{a\in A} q^a.
\]
Then
\[
k = f_A(1) = |A|
\]
and
\[
\frac{f_A(q)^2 + f_A(q^2)}{2} = \sum_{j\in 2A}r_{2,A}(j)q^j.
\]
If $[0,n-1] \subseteq 2A$, then $r_{2,A}(j) \geq 1$ for all $0 \leq j\leq n-1$.
Hence there exist integers $\delta(j) \geq 0$ such that
\[
\frac{f_A(q)^2 + f_A(q^2)}{2} = 1 + q + q^2 + \cdots + q^{n-1} +
\sum_{j\in 2A}\delta(j)q^j,
\]
where
\[
\delta(j) = \left\{
\ba{ll}
r_{2,A}(j) - 1 & \text{ if $j \in \{0,1,\ldots, n-1\}$,}\\
r_{2,A}(j) & \text{otherwise.}
\ea
\right.
\]
Let
\[
\Delta(q) = \sum_{j\in 2A}\delta(j)q^j.
\]
Then $\Delta(q) \geq 0$ for $q \geq 0,$ and
\beq  \label{rohr:moseriden}
\frac{f_A(q)^2 + f_A(q^2)}{2} = 1 + q + q^2 + \cdots + q^{n-1} + \Delta(q).
\eeq
Evaluating the generating function identity~(\ref{rohr:moseriden}) at
$q = 1$, we obtain
\beq       \label{rohr:main}
\frac{k^2+k}{2} = n + \Delta(1).
\eeq
Since $\Delta(1) \geq 0,$ we have
\[
n \leq \frac{k^2}{2} + O(k).
\]
The strategy is to find a lower bound for $\Delta(1)$ of the form
\[
\Delta(1) \geq ck^2 + O(k)
\]
for some $c > 0,$ and deduce
\[
n \leq \left(\frac{1}{2} - c\right)k^2 + O(k).
\]

We obtain a simple combinatorial lower bound for $\Delta(1)$ by noting that
if $a_1, a_2 \in A$ and $n/2 \leq a_1 \leq a_2$, then $a_1+a_2 \geq n.$
Let $\ell$ denote the number of elements $a \in A$ such that $a \geq n/2.$
Then
\beq   \label{rohr:combineq}
\Delta(1) \geq \sum_{j\geq n}\delta(j)
= \sum_{j \geq n} r_{2,A}(j) \geq \frac{\ell(\ell+1)}{2} \geq \frac{\ell^2}{2}.
\eeq
Let
\[
\omega = e^{2\pi i/n}
\]
be a primitive $n$th root of unity.
Let $r$ be an integer not divisible by $n$.  Then
\[
1 + \omega^r + \omega^{2r} + \cdots + \omega^{(n-1)r} = 0
\]
and so
\[
\frac{ f_A(\omega^r)^2 + f_A(\omega^{2r}) }{2}
= 1 + \omega^r + \omega^{2r} + \cdots + \omega^{(n-1)r} + \sum_{j}\delta(j)\omega^{jr}  = \Delta(\omega^r).
\]
Applying the triangle inequality, we obtain
\[
\Delta(1) \geq |\Delta(\omega^r)| = \frac{ \left|f_A(\omega^r)^2 +
f_A(\omega^{2r})\right| }{2} \geq \frac{ |f_A(\omega^r)|^2 - k }{2}.
\]
Let
\[
M = \max\{ |f_A(\omega^r)| : r \not\equiv 0 \pmod{n}\}.
\]
Then
\beq    \label{rohr:Mkineq}
0 \leq M \leq k
\eeq
and
\beq  \label{rohr:Mineq}
\Delta(1) \geq \frac{ M^2 - k }{2}.
\eeq

We can also obtain an analytic lower bound for $\Delta(1).$
For all integers $r$ not divisible by $n$, we have
\[
M \geq |f_A(\omega^r)| = \left| \sum_{a\in A}e^{2\pi i ra/n}  \right|
= \left| \sum_{a\in A} \cos(2\pi ra/n) + i \sin(2\pi ra/n)\right|,
\]
and so
\[
\left| \sum_{a\in A} \cos(2\pi ra/n) \right| \leq M
\]
and
\[
\left| \sum_{a\in A} \sin(2\pi ra/n)\right| \leq M.
\]

Let $\varphi(t)$ be a function with period 1 and with a Fourier series
\[
\varphi(t) = \sum_{r=0}^{\infty} a_r \cos (2\pi rt) +  \sum_{r=1}^{\infty}b_r \sin (2\pi rt)
\]
whose Fourier coefficients converge absolutely, that is,
\[
\sum_{r=0}^{\infty} |a_r| + \sum_{r=1}^{\infty} |b_r| < \infty.
\]

Let
\[
C = \sum_{\genfrac{}{}{0pt}{}{r=0}{n \mid r}}^{\infty} |a_r|.
\]
For any integer $a$ we have
\[
\begin{split}
\sum_{a\in A} \varphi\left(\frac{a}{n}\right)
& = \sum_{a\in A} \sum_{r=0}^{\infty} a_r \cos (2\pi ra/n) +
 \sum_{a\in A} \sum_{r=0}^{\infty}b_r \sin (2\pi ra/n) \\
& = \sum_{r=0}^{\infty} a_r \sum_{a\in A} \cos (2\pi ra/n) +
\sum_{r=1}^{\infty} b_r \sum_{a\in A}\sin (2\pi ra/n) \\
& = \sum_{\genfrac{}{}{0pt}{}{r=0}{n \nmid r}}^{\infty} a_r\sum_{a\in A}
\cos (2\pi ra/n) +
\sum_{\genfrac{}{}{0pt}{}{r=1}{n \nmid r}}^{\infty} b_r \sum_{a\in A}
\sin (2\pi ra/n) + k \sum_{\genfrac{}{}{0pt}{}{r=0}{n \mid r}}^{\infty} a_r,\\
\end{split}
\]
and so
\[
\left| \sum_{a\in A} \varphi\left(\frac{a}{n}\right) \right| \leq
M \sum_{\genfrac{}{}{0pt}{}{r=0}{n \nmid r}}^{\infty} (|a_r| + |b_r|) + k C.
\]

Let $\alpha_1$ and $\alpha_2$ be real numbers such that
\[
\varphi(t) \geq \alpha_1 \qquad\text{for $0 \leq t <1/2$}
\]
and
\[
\varphi(t) \geq \alpha_2\qquad\text{for $1/2 \leq t < 1.$}
\]
Recall that $\ell$ denotes the number of elements $a \in A$ such that
$n/2 \leq a \leq n-1.$  Then
\[
\sum_{a\in A} \varphi\left(\frac{a}{n}\right) \geq (k-\ell)\alpha_1 + \ell\alpha_2 = k\alpha_1 - (\alpha_1 - \alpha_2)\ell.
\]
We obtain the inequality
\beq      \label{rohr:FSineq}
k\alpha_1 - (\alpha_1 - \alpha_2)\ell \leq M
\sum_{\genfrac{}{}{0pt}{}{r=0}{n \nmid r}}^{\infty} (|a_r| + |b_r|) + k C.
\eeq
In this way, the function $\varphi(t)$ produces a lower bound for $M$, which, by~(\ref{rohr:Mineq}), gives a lower bound for $\Delta(1).$

Moser applied inequality~(\ref{rohr:FSineq}) to the function
\[
\varphi(t) = \frac{1}{2}\cos(4\pi t) + \sin(2\pi t),
\]
whose nonzero Fourier coefficients are $a_2 = 1/2$ and $b_1 = 1.$
Then $C = 0$ for $n \geq 3,$ and
\[
\left| \sum_{a\in A} \varphi\left(\frac{a}{n}\right) \right| \leq \frac{3M}{2}.
\]
The function $\varphi(t)$ satisfies the inequality
\[
\varphi(t) \geq \left\{
\ba{ll} \frac{1}{2} & \text{for $0 \leq t < 1/2$}\\
-\frac{3}{2} & \text{for $1/2 \geq t < 1$,}
\ea
\right.
\]
and so
\[
\sum_{a\in A} \varphi\left(\frac{a}{n}\right)
\geq \frac{k-\ell}{2} - \frac{3\ell}{2} = \frac{k-4\ell}{2}.
\]
This implies that
\[
M \geq \frac{2}{3}\left| \sum_{a\in A} \varphi\left(\frac{a}{n}\right) \right| \geq
\frac{k-4\ell}{3},
\]
and we obtain the analytic lower bound
\[
\Delta(1) \geq \frac{(k-4\ell)^2}{18} - \frac{k}{2}.
\]
Recalling the combinatorial lower bound~(\ref{rohr:combineq})
\[
\Delta(1) \geq \frac{M^2-k}{2},
\]
we obtain
\[
\Delta(1)
\geq \max\left\{ \frac{(k-4\ell)^2}{18}, \frac{\ell^2}{2} \right\} - \frac{k}{2} = \frac{k^2}{98} - \frac{k}{2}.
\]
Inserting this into inequality~(\ref{rohr:main}), we obtain
\[
\frac{k^2+k}{2} = n + \Delta(1) \geq n + \frac{k^2}{98} - \frac{k}{2},
\]
and so
\[
n \leq  \left(\frac{1}{2} - \frac{1}{98}\right) k^2 + k \leq 0.4898 k^2+k.
\]

\section{Fourier series in two variables}

We shall modify Moser's method to obtain a better lower bound for $\Delta(1).$
We use the same notation as in the previous section.
In particular, $\ell$ denotes the number of integers $a \in A$ such that
$a \geq n/2.$  Let $L$ denote the number of pairs $(a_1,a_2) \in A \times A$
such that $a_1+a_2 \geq n.$
Then $L \geq \ell^2$, and $k^2 - L$ is the number of pairs
$(a_1,a_2) \in A \times A$ such that $a_1+a_2 \leq n-1$.
We have the combinatorial lower bound
\beq    \label{rohr:Lineq}
\Delta(1) \geq \sum_{j\geq n} r_{2,A}(n) = \frac{L+\ell}{2} \geq \frac{L}{2}.
\eeq

Let $\varphi(t_1,t_2)$ be a function with period 1 in each variable and with a Fourier series
\[
\varphi(t_1,t_2)
= \sum_{r_1\in \Z}\sum_{r_2\in \Z} \hat{\varphi}(r_1,r_2)e^{2\pi i r_1t_1}e^{2\pi i r_2t_2}
\]
whose Fourier coefficients converge absolutely, that is,
\[
\sum_{r_1\in \Z} \sum_{r_2\in \Z} \left| \hat{\varphi}(r_1,r_2) \right| < \infty.
\]
We choose $\varphi(t_1,t_2)$ with zero mean, that is,
\[
\hat{\varphi}(0,0) = \int_0^1\int_0^1 \varphi(t_1,t_2) dt_1 dt_2 = 0.
\]
Let
\[
R_1 = \{ (t_1,t_2) \in [0,1)\times [0,1) : t_1+t_2 < 1 \}
\]
and let
\[
R_2 = \{ (t_1,t_2) \in [0,1)\times [0,1) : t_1+t_2 \geq 1 \}
\]
If $a_1,a_2 \in A$ and $a_1+a_2 \leq n-1$, then $(a_1/n, a_2/n) \in R_1$.
If $a_1+a_2 \geq n$, then $(a_1/n, a_2/n) \in R_2.$

Let $\alpha_1$ and $\alpha_2$ be real numbers such that
\[
\varphi(t_1,t_2) \geq \alpha_1 \qquad\text{for $(t_1,t_2) \in R_1$}
\]
and
\[
\varphi(t_1,t_2) \geq \alpha_2 \qquad\text{for $(t_1,t_2) \in R_2.$}
\]
We choose the function $\varphi(t_1,t_2)$ such that
\[
\alpha_1 > \alpha_2.
\]
Then
\beq   \label{rohr:A}
\sum_{a_1\in A} \sum_{a_2\in A}\varphi\left(\frac{a_1}{n},\frac{a_2}{n}\right)
\geq (k^2-L)\alpha_1 + L\alpha_2 = \alpha_1 k^2 - (\alpha_1 - \alpha_2)L.
\eeq
We can rewrite this sum as follows:
\[
\begin{split}
\sum_{a_1\in A} \sum_{a_2\in A}\varphi\left(\frac{a_1}{n},\frac{a_2}{n}\right)
& = \sum_{a_1\in A} \sum_{a_2\in A} \sum_{r_1\in \Z} \sum_{r_2\in \Z}
\hat{\varphi}(r_1,r_2)e^{2\pi i r_1 a_1/n}e^{2\pi i r_2 a_2/n} \\
& = \sum_{r_1\in \Z}\sum_{r_2\in \Z} \hat{\varphi}(r_1,r_2)
\sum_{a_1\in A}e^{2\pi i r_1 a_1/n} \sum_{a_2\in A}e^{2\pi i r_2 a_2/n} \\
& = \sum_{r_1\in \Z}\sum_{r_2\in \Z} \hat{\varphi}(r_1,r_2)
f_A(\omega^{r_1})f_A(\omega^{r_2}).
\end{split}
\]
Consider the partition of the integer lattice $\Z^2 = S_0 \cup S_1 \cup S_2$:
\[
\begin{split}
S_0 = & \{ (r_1,r_2) \in \Z^2 : r_1 \equiv r_2 \equiv 0 \pmod{n} \}  \\
S_1  = & \{ (r_1,r_2) \in \Z^2 : r_1 \equiv 0\pmod{n}, r_2 \not\equiv 0 \pmod{n} \}  \\
 & \cup \{ (r_1,r_2) \in \Z^2 : r_1 \not\equiv 0\pmod{n}, r_2 \equiv 0 \pmod{n} \}  \\
S_2 = & \{ (r_1,r_2) \in \Z^2 : r_1 \not\equiv 0 \pmod{n},r_2 \not\equiv 0 \pmod{n} \}.
\end{split}
\]
We define $C_0, C_1,$ and $C_2$ by
\[
C_i = \sum_{(r_1,r_2) \in S_i} \left| \hat{\varphi}(r_1,r_2) \right|.
\]
Recall that $|f_A\left( \omega^r \right)| \leq M$ if $r$ is not divisible by $n$
and $|f_A\left( \omega^r \right)| \leq k$ if $r$ is divisible by $n$.
Then
\beq           \label{rohr:B}
\left|  \sum_{a_1\in A} \sum_{a_2\in A}
\varphi\left(\frac{a_1}{n},\frac{a_2}{n}\right) \right|
\leq C_0 k^2 + C_1kM + C_2 M^2.
\eeq
Combining inequalities~(\ref{rohr:A}) and~(\ref{rohr:B}), we obtain
\[
\alpha_1 k^2 - (\alpha_1 - \alpha_2)L \leq C_0 k^2 + C_1kM + C_2 M^2.
\]
Since $\alpha_1 > \alpha_2,$ we have
\[
L \geq \frac{(\alpha_1 - C_0 )k^2 - C_1kM - C_2 M^2}{\alpha_1 - \alpha_2}.
\]
We define
\[
\mu = \frac{M}{k}.
\]
Since $0 \leq M \leq k,$ we have
\[
0 \leq \mu \leq 1.
\]
By inequality~(\ref{rohr:Lineq}), we have $2\Delta(1) \geq L,$ and so
\[
\frac{2\Delta(1)}{k^2} \geq \frac{L}{k^2} \geq \frac{(\alpha_1 - C_0 ) - C_1\mu - C_2 \mu^2}{\alpha_1 - \alpha_2}.
\]
By inequality~(\ref{rohr:Mineq}), we also have $2\Delta(1) \geq M^2 - k$, and so
\beq      \label{rohr:C}
\frac{2\Delta(1)}{k^2} \geq \max\left( \mu^2,\frac{(\alpha_1 - C_0 ) - C_1\mu - C_2 \mu^2}{\alpha_1 - \alpha_2}\right) - \frac{1}{k}.
\eeq

Since the series of Fourier coefficients of $\varphi(t_1,t_2)$ converges absolutely
and since $\hat{\varphi}(0,0) = 0,$ we can arrange the Fourier series in the form of a sum over concentric squares
\[
\sum_{R=1}^{\infty} \sum_{\max(|r_1|,|r_2|) = R} \hat{\varphi}(r_1,r_2)
e^{2\pi i r_1 t_1}e^{2\pi i r_2 t_2}.
\]
For any $\varepsilon > 0$ there exists an integer $N = N(\varepsilon)$
such that
\[
\sum_{n=N}^{\infty} \sum_{\max(|r_1|,|r_2|) = n}
\left|\hat{\varphi}(r_1,r_2)\right| < \epsilon(\alpha_1 - \alpha_2).
\]
For all $n \geq N,$ we shall approximate the sums $C_0, C_1,$ and $C_2$ by 0, $C_{\text{axial}},$ and $C_{\text{main}},$ respectively, where
\[
C_{\text{axial}} = \sum_{r \in \Z \atop r\neq 0}
\left( |\hat{\varphi}(0,r)| +  |\hat{\varphi}(r,0)| \right)
\]
and
\[
C_{\text{main}} = \sum_{r_1 \in \Z \atop r_1\neq 0}
\sum_{r_2 \in \Z \atop r_2\neq 0} |\hat{\varphi}(r_1,r_2)|.
\]
Then
\[
\begin{split}
\left|
(\alpha_1 - C_0 ) - C_1\mu - C_2 \mu^2  \right.
& \left.
- \left(
\alpha_1 - C_{\text{axial}}\mu - C_{\text{main}} \mu^2
\right)
\right|  \\
& = \left| C_0 + (C_1 - C_{\text{axial}}) \mu
+ (C_2 - C_{\text{main}})\mu^2\right| \\
& \leq |C_0| + |C_1 - C_{\text{axial}}| + |C_2 - C_{\text{main}}| \\
& \leq \sum_{\max(|r_1|,|r_2|) \geq N} \left| \hat{\varphi}(r_1,r_2)\right| \\
& < \epsilon(\alpha_1 - \alpha_2),
\end{split}
\]
and so
\[
\left|
\left(
\frac{(\alpha_1 - C_0 ) - C_1\mu  - C_2 \mu^2}{\alpha_1 - \alpha_2} \right)
- \left(  \frac{\alpha_1  - C_{\text{axial}}\mu - C_{\text{main}} \mu^2}{\alpha_1 - \alpha_2}
\right)
\right| < \varepsilon.
\]
It follows from inequality~(\ref{rohr:C}) that
\[
\frac{2\Delta(1)}{k^2} \geq \max\left( \mu^2,\frac{\alpha_1 - C_{\text{axial}}\mu - C_{\text{main}} \mu^2}{\alpha_1 - \alpha_2}\right) - \varepsilon - \frac{1}{k}.
\]
Let
\beq             \label{rohr:rho}
\rho = \inf_{0 \leq \mu \leq 1} \max\left( \mu^2,  \frac{\alpha_1  - C_{\text{axial}}\mu - C_{\text{main}} \mu^2}{\alpha_1 - \alpha_2} \right).
\eeq

From (\ref{rohr:C}) and the definition of $\rho$ in (\ref{rohr:rho}), we
now have
\[
\frac{2\Delta(1)}{k^2} \geq \rho - \varepsilon - \frac{1}{k}.
\]
Applying identity~(\ref{rohr:main}), we obtain
\[
\frac{k^2+k}{2}  = n + \Delta(1) \geq n + \frac{(\rho-\varepsilon)k^2-k}{2}.
\]
Therefore,
\[
n \leq \left( \frac{1-\rho + \varepsilon}{2}\right)k^2 + k,
\]
where the number $\rho$ depends only on the function $\varphi(t_1,t_2)$
and $\epsilon >0$ can be arbitrary small.

It is clear that we always have $\rho \geq 0$, and that
$\rho > 0$ if and only if $\alpha_1 > 0$. It is also clear
that when $\alpha_1 \geq 0$, we have
$\rho = \xi^2$, where $\xi$ is the unique solution
in $[0,1]$ to the quadratic equation
$$ \xi^2 = \frac{\alpha_1  - C_{\text{axial}}\xi - C_{\text{main}} \xi^2}
{\alpha_1 - \alpha_2},
$$
i.e.,
$$ (\alpha_1-\alpha_2+C_{\text{main}})\xi^2 + C_{\text{axial}}\xi - \alpha_1
= 0,
$$
which yields the formula
\begin{equation}
\rho =
\left( \frac{-C_{\text{axial}}+
\sqrt{C^2_{\text{axial}}+4\alpha_1(\alpha_1-\alpha_2+C_{\text{main}})}}
{2\alpha_1(\alpha_1-\alpha_2+C_{\text{main}})}\right)^{\!2}.
\end{equation}

Hence we have an optimization problem in which we maximize
$\rho$ over all real valued functions $\varphi$ defined on the
unit square $[0,1)^2$ such that $\varphi$ has zero mean and
$\varphi > 0$ on $R_1$.
We do not know the optimal function for this problem, but we have found
a simple piecewise polynomial function that improves Klotz's upper
bound for $n(2,k)$. Before we proceed to the main result of this paper,
which also includes the definition of this function,
let us present some of the heuristics which have lead us to our ``educated
guess.''

First, without loss of generality, we may assume that
$\alpha_1 = 1$. Note that we then necessarily have
$$ \frac{1}{2}\alpha_1 + \frac{1}{2}\alpha_2 \leq
\iint_{R_1} \varphi(t_1,t_2) dt_1 dt_2 +
\iint_{R_2} \varphi(t_1,t_2) dt_1 dt_2 = 0 $$
so that $\alpha_2 \leq -1$. We also have
$$ C_{\text{axial}} \geq
\left | \sum_{r} \left (\hat \varphi(r,0) + \hat \varphi(0,r) \right ) \right |
= \left | \int_0^1 \left (\varphi(0,t) + \varphi(t,0) \right ) dt \right |
\geq 2,$$
and
$$C_{\text{main}} \geq \left | \sum_{r}
\hat \varphi(r,r)  \right | =
\left | \int_0^1 \varphi(t,1-t) dt \right | \geq 1.$$
In any case we are interested in the positive root $\xi_{\kappa,\tau}$ of
the equation
$$ \kappa \xi^2 + \tau \xi - 1 = 0$$
where $\kappa = (1-\alpha_2 + C_{\text{main}}) \geq 3$
and $\tau = C_{\text{axial}} \geq 2$. Clearly, the smaller $\kappa$
and $\tau$ are, the larger this root will be. The bounds
$\kappa \geq 3$ and $\tau \geq 2$ already
imply that $\xi_{\kappa,\tau} \leq \frac{1}{3}$, hence
$\rho = \xi^2_{\kappa,\tau} \leq \frac{1}{9}$.
In reality, $\alpha_2 < -1$ because equality can happen only if
$\varphi$ is constant on both $R_1$ and $R_2$, in which case
$\hat \varphi$ is not absolutely summable. This results
in the heuristic that if we try to push $\alpha_2$ close to $-1$, then
$C_{\text{axial}}$ and $C_{\text{main}}$ will become large,
and conversely if we try to push $C_{\text{axial}}$ and
$C_{\text{main}}$ close to their respective minimum values, then
$\varphi$ may not be bounded from below on $R_2$ by a small value.
The right trade-off between these two competing quantities will
result in the solution of this optimization problem.

It is interesting to note that the value of $\rho$ is fairly robust
with respect to variations in $\kappa$ and $\tau$, which we will only
be able to estimate but not compute exactly. The following lemma
gives an explicit estimate for this purpose:

\begin{lemma}\label{root-variation}
Let $\xi_{\kappa,\tau}$ and $\xi_{\kappa_0,\tau_0}$ be the respective
positive roots
of the equations $ \kappa \xi^2 + \tau \xi - 1 = 0$ and
$ \kappa_0 \xi^2 + \tau_0 \xi - 1 = 0$. Let
$\rho = \xi_{\kappa,\tau}^2$ and $\rho_0 = \xi_{\kappa_0,\tau_0}^2$.
If $\min(\kappa,\kappa_0) \geq 3$ and $\min(\tau, \tau_0) \geq 2$, then
\begin{equation}
|\rho - \rho_0| \leq \frac{1}{54}|\kappa-\kappa_0|+\frac{1}{18}|\tau-\tau_0|.
\end{equation}
\end{lemma}

The proof of this lemma is given in the Appendix.
Now we can state and prove the main theorem of this paper.

\bt
\[
\limsup_{n\rightarrow\infty}\frac{n(2,k)}{k^2} \leq 0.4789.
\]
\et

\begin{proof}
We define the function $\varphi(t_1,t_2)$ on the unit square $[0,1)^2$
by
\begin{equation}\label{def-varphi}
\varphi(t_1,t_2) = \left\{
\ba{ll}
1, & (t_1,t_2) \in R_1  \\
1-40(1-t_1)(1-t_2)\left(1-(2-t_1-t_2)^6\right), &(t_1,t_2) \in R_2.
\ea
\right.
\end{equation}

\begin{figure}[t]
\begin{center}
\includegraphics[height=3.5in]{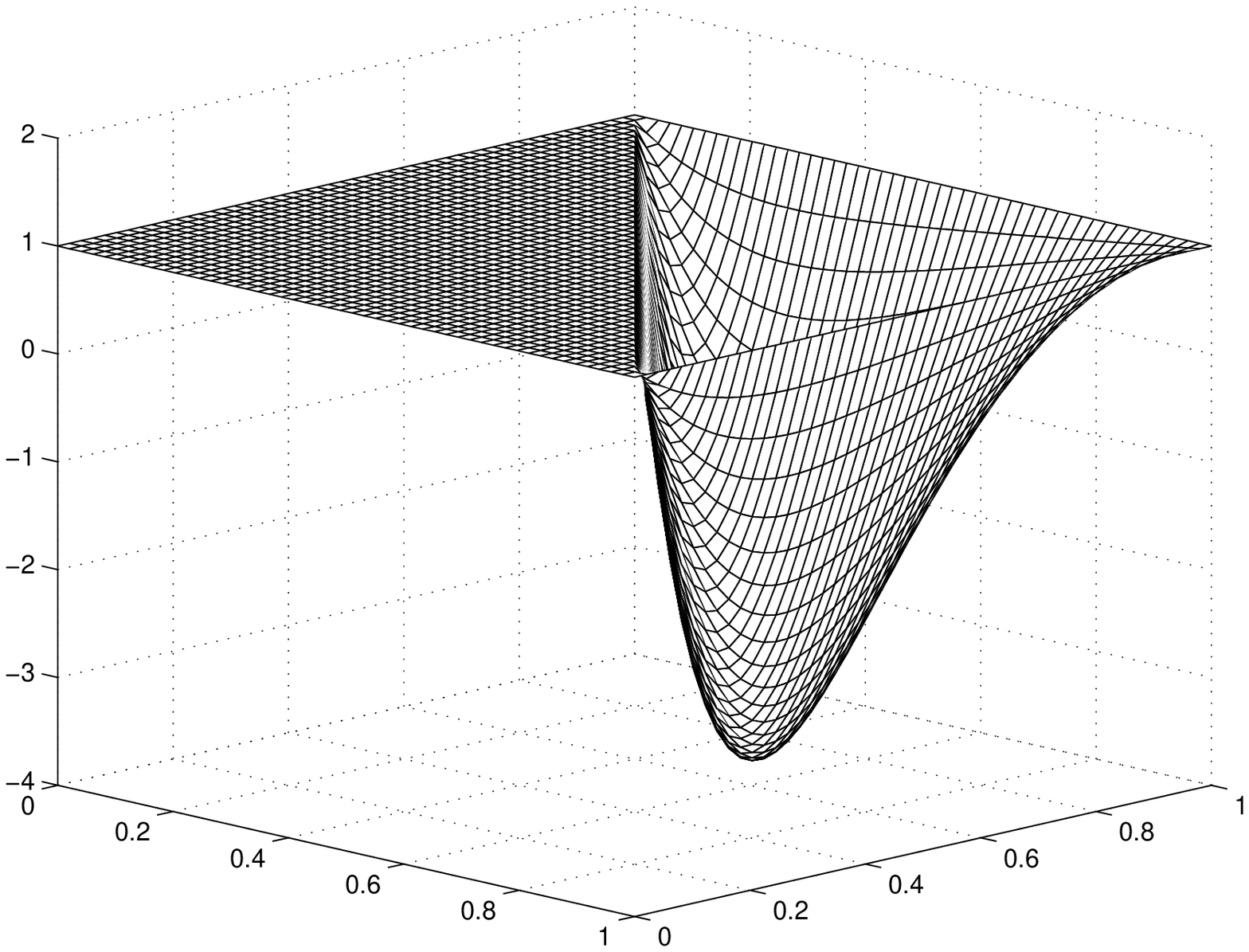}
\end{center}
\caption{\label{graph} }
\end{figure}

Then
\[
\alpha_1 = 1.
\]
Computation of the three other parameters used in formula~(\ref{rohr:rho})
for $\rho$ yields
\[
\alpha_2 = 1 - \frac{15}{2^{5/3}} = -3.72470\ldots,
\]
\[
2.90278 \leq C_{\text{axial}} \leq 2.90289,
\]
and
\[
4.75145 \leq C_{\text{main}} \leq 4.76146.
\]
Taking $\kappa_0 = 1+3.72471+4.76146 = 9.48617$, and
$\tau_0 = 2.90289$, we obtain $\rho_0 > 0.04240$,
$|\kappa - \kappa_0| < 0.01002$ and $|\tau - \tau_0| < 0.00011$, so that
$|\rho - \rho_0| < 0.0002$. Hence
\[
\rho \geq  \rho_0 - |\rho-\rho_0| > 0.0422,
\]
and consequently choosing $\epsilon$ sufficiently small,
\[
n \leq 0.4789\, k^2+k.
\]
The details of the computations are in the Appendix to this paper.
This completes the proof.
\end{proof}

\section{Open problems}
A major open problem concerning the extremal function
\[
n(2,k) = \max\{ n(2,A) : A \subseteq \N_0 \text{ and } |A| \leq k\}.
\]
is to compute $\liminf_{n\rightarrow\infty}n(2,k)/k^2$ and $\limsup_{n\rightarrow\infty} n(2,k)/k^2,$ and to determine if the limit
\[
\lim_{n\rightarrow\infty}\frac{n(2,k)}{k^2}
\]
exists.  We have no conjecture about the existence of this limit, nor about the values of the $\liminf$ and $\limsup.$

It is also difficult to compute the exact values of the function $n(2,k).$

We can generalize the extremal functions $n(2,A)$ and $n(2,k)$ as follows.
Let $A$ be a finite set of integers, and let $m(2,A)$ denote the largest integer $n$ such that the sumset $2A$ contains $n$ consecutive integers.  Let
\[
m(2,k) = \max\{ m(2,A) : A \subseteq \Z \text{ and } |A| \leq k\}.
\]
Let $\ell(2,A)$ denote the largest integer $n$ such that the sumset $2A$ contains an arithmetic progression of length $n$, and let
\[
\ell(2,k) = \max\{ \ell(2,A) : A \subseteq \Z \text{ and } |A| \leq k\}.
\]
We can also define the extremal function
\[
n'(2,k) = \max\{ n(2,A) : A \subseteq \Z \text{ and } |A| \leq k\}.
\]
Then
\[
n(2,A) \leq n'(2,A) \leq m(2,A) \leq \ell(2,A),
\]
and so
\[
n(2,k) \leq n'(2,k) \leq m(2,k) \leq \ell(2,k).
\]
For any integer $t$ and set $A$, we have the translation
$A+t = \{a+t:a\in A\}.$
The functions $\ell$ and $m$ are translation invariant, that is,
$\ell(2,A+t) = \ell(2,A)$ and  $m(2,A+t) = m(2,A)$.
We also have the trivial upper bound $\ell(2,k) \leq \binom{k+1}{2}$,
but it is an open problem to obtain nontrivial upper bounds for any of
the extremal functions $n'(2,k)$, $m(2,k)$, or $\ell(2,k)$.

\section*{Appendix}
We describe here the computations.

\subsection*{Proof of Lemma \ref{root-variation}}
We start with the formula
$$\xi_{\kappa,\tau} = \frac{-\tau+\sqrt{\tau^2+4\kappa}}{2\kappa}
 = \frac{2}{\tau+\sqrt{\tau^2+4\kappa}}.$$
We next evaluate the partial derivatives of $\xi_{\kappa,\tau}$
with respect to $\kappa$ and $\tau$:
\begin{eqnarray}
\frac{\partial \xi_{\kappa,\tau}}{\partial \kappa}
& = & -\frac{4}{\sqrt{\tau^2+4\kappa}\,(\tau + \sqrt{\tau^2+4\kappa})^2}, \\
\frac{\partial \xi_{\kappa,\tau}}{\partial \tau}
& = & -\frac{2}{\sqrt{\tau^2+4\kappa}\,(\tau + \sqrt{\tau^2+4\kappa})}
\end{eqnarray}
from which it follows that in the set $\{(\kappa,\tau) :
\kappa \geq 3, \; \tau \geq 2 \}$, we have
$\left |\frac{\partial \xi_{\kappa,\tau}}{\partial \kappa} \right |
\leq \frac{1}{36}$ and
$\left |\frac{\partial \xi_{\kappa,\tau}}{\partial \tau} \right |
\leq \frac{1}{12}$. These bounds then imply
$$
\left |\xi_{\kappa,\tau} - \xi_{\kappa_0,\tau_0} \right |
\leq \frac{1}{36} |\kappa - \kappa_0| +  \frac{1}{12} |\tau - \tau_0|.
$$
We then note that
$\xi_{\kappa,\tau}  \leq \frac{1}{3}$, which yields
$$
|\rho - \rho_0|  = \left |\xi_{\kappa,\tau} - \xi_{\kappa_0,\tau_0} \right |
\left |\xi_{\kappa,\tau} +  \xi_{\kappa_0,\tau_0} \right |
\leq \frac{2}{3}  \left |\xi_{\kappa,\tau} - \xi_{\kappa_0,\tau_0} \right |,
$$
hence the result of the lemma.
\hfill $\Box$

\subsection*{Absolute summability of $\hat \varphi$}
While the result explained in this subsection is elementary, we will
provide a certain amount of detail in its derivation because our main
concern is more than  absolute summability of $\hat \varphi$. We would
like to provide explicit estimates on the rate of the convergence;
this will be necessary later in the section when we will analyze
the accuracy of the numerical computation of the
constants $C_{\text{main}}$ and $C_{\text{axial}}$.

\begin{lemma}\label{fourier-decay-1}
Let $f$ be a smooth function on $[0,1]$. Then for all $L \geq 0$
and $n \not= 0$, the following formula holds:
\begin{equation}
\hat f(n) = \sum_{k=0}^{L} \frac{f^{(k)}(0)-f^{(k)}(1)}{(2\pi i n)^{k+1}}
+ \frac{\widehat{f^{(L+1)}}(n)}{(2\pi i n)^{L+1}}.
\end{equation}
\end{lemma}

\begin{proof}
The case $L=0$ follows from integration by parts and
the general case follows from iterating this result.
\end{proof}

\begin{theorem}
Let $F$ be a smooth function of two real
variables and assume that $F$ vanishes on the boundary of $R_2$, i.e.,
$$
F(t,1-t) = F(1,t) = F(t,1) = 0,
~~\mbox{ for all } t \in [0,1].$$
Define
$$
\Psi_F(t_1,t_2) = \left\{
\ba{ll}
0, & \text{if $(t_1,t_2) \in R_1$,}\\
F(t_1,t_2), & \text{if $(t_1,t_2) \in R_2$.}
\ea
\right.
$$
Then the Fourier series expansion of
$\Psi_F$ is absolutely convergent.
\end{theorem}

\par \noindent {\bf Note:} Later we will simply set $\varphi = \Psi_F + 1$.

\begin{proof}
We will prove this result by
deriving a suitable decay estimate on $|\hat \Psi_F(r_1,r_2)|$, where
$$
\hat \Psi_F(r_1,r_2) = \int_0^1 e^{-2\pi i r_1 t_1} \left \{
\int_{1-t_1}^1  F(t_1,t_2) e^{-2\pi i r_2 t_2}\, dt_2 \right\} dt_1.
$$

\subsubsection*{The case $r_1 = 0$ or $r_2 =0$.}
Due to the symmetry on the assumptions on $F$,
it suffices to consider only one of these cases. Let us assume that
$r_2 = 0$. Define
$$ J_1(t_1) = \int_{1-t_1}^1  F(t_1,t_2)\, dt_2 $$
so that
\begin{equation} \label{on-the-axes}
\hat \Psi_F(r_1,0) = \hat J_1(r_1).
\end{equation}
Clearly we have $J_1(0) = J_1(1) = 0$.
Setting $L=1$ and $f = J_1$
in Lemma \ref{fourier-decay-1}, we see that for $r_1 \not= 0$
\begin{equation}\label{axial-decay-1}
|\hat \Psi_F(r_1,0)| \leq
\frac{1}{|2 \pi r_1|^2}\left( |J_1'(0) - J_1'(1)| + \int_0^1 |J_1''| \right)
= O\left (\frac{1}{r_1^2} \right ).
\end{equation}
With a similar estimate for $ |\hat \Psi_F(0,r_2)|$, we have
$$\sum_{r \not= 0}
\left( |\hat \Psi_F(r,0)|+ |\hat \Psi_F(0,r)| \right) < \infty.$$

\subsubsection*{The case $r_1 \not= 0$ and $r_2 \not =0$.}
We will derive a general formula for $\hat \Psi_F(r_1,r_2)$. To do this,
we momentarily forget that $F$ vanishes on the boundary of $R_2$, and
for $t \in [0,1]$, define the following functions:
$$
\begin{array}{rcrcrcr}
g_{0}(t) &=& F(t,1-t), & &
h_{0}(t) &=& F(t,1),  \\
g_{1}(t) &=& (\partial_2 F)(t,1-t),  & &
h_{1}(t) &=& (\partial_2 F)(1,t), \\
g_{2}(t) &=& (\partial_1 \partial_2 F)(t,1-t), & &
h_{2}(t) &=& (\partial_1 \partial_2 F)(t,1), \\
g_{3}(t) &=& (\partial_1 \partial^2_2 F)(t,1-t), & &
h_{3}(t) &=& (\partial_1 \partial^2_2 F)(1,t).
\end{array}
$$
We start with the formula for $\hat \Psi_F(r_1,r_2)$ above.
Integrating by parts in the second variable, we obtain
\begin{eqnarray}
\hat \Psi_F(r_1,r_2)
& = &
\int_0^1 dt_1\, e^{-2\pi i r_1 t_1} \bigg \{
\left [\frac{e^{2 \pi i r_2 t_2}}{-2 \pi i r_2} F(t_1,t_2)
\right ]_{t_2=1-t_1}^{t_2 =1} \nonumber \\
& & \qquad \qquad \qquad \qquad
-  \int_{1-t_1}^1 dt_2\, \frac{e^{-2\pi i r_2 t_2}}{-2 \pi i r_2}
(\partial_2 F)(t_1,t_2) \bigg \}\nonumber \\
& = &
\frac{1}{(2 \pi i) r_2}
\Big(\hat g_0(r_1-r_2) -  \hat h_0(r_1) +
\hat \Psi_{\partial_2 F}(r_1,r_2) \Big).  \nonumber
\end{eqnarray}
We apply the same method to
$\hat \Psi_{\partial_2 F}(r_1,r_2)$, but integrate by parts in the
first variable. This results in
\begin{equation*}
\begin{split}
\hat \Psi_F(r_1,r_2)
& = \frac{1}{(2 \pi i) r_2}
\Big(\hat g_0(r_1-r_2) -  \hat h_0(r_1) \Big) \\
& \quad + \frac{1}{(2 \pi i)^2 r_1 r_2}
\Big(\hat g_1(r_1-r_2) - \hat h_1(r_2) +
\hat \Psi_{\partial_1 \partial_2 F}(r_1,r_2) \Big)
\end{split}
\end{equation*}
We repeat the first two steps in the same order, which gives us
\begin{equation}
\begin{split}
\hat \Psi_F(r_1,r_2)
& = \frac{1}{(2 \pi i) r_2}
\Big(\hat g_0(r_1-r_2) -  \hat h_0(r_1) \Big) \\
&\quad + \frac{1}{(2 \pi i)^2 r_1 r_2}
\Big(\hat g_1(r_1-r_2) - \hat h_1(r_2) \Big) \\
&\quad + \frac{1}{(2 \pi i)^3 r_1 r^2_2}
\Big(\hat g_2(r_1-r_2) -  \hat h_2(r_1) \Big) \\
&\quad + \frac{1}{(2 \pi i )^4 r^2_1 r^2_2}
\Big(\hat g_3(r_1-r_2) -  \hat h_3(r_2)
+ \hat \Psi_{\partial^2_1 \partial^2_2 F}(r_1,r_2) \Big)
\end{split}
\end{equation}

Note that from our assumptions on $F$, we have $g_0 = h_0 = h_1 = 0$.
We will have two subcases:
\begin{enumerate}
\item $r_1 = r_2 = r$.
In this case, we easily see from the second formula above that
\begin{equation}\label{diagonal-decay-1}
|\hat \Psi_F(r,r)| \leq
\frac{|\hat g_1(0)| + \| \Psi_{\partial_1 \partial_2 F} \|_1}{4\pi^2 r^2}
\end{equation}
\item $r_1 \not= r_2$.
This case is slightly more subtle. We first note that $g_1(1) = h_1(0) = 0$.
It is also true that $g_1(0) = (\partial_2 F)(0,1)$. To see this, note that
$(\partial_1 F)(0,1) = 0$ and $\nabla F(0,1) \cdot (1,-1) = 0$,
both of which follow from the fact that $F$ vanishes on the boundary of
$R_2$. The function $g_1$ being smooth otherwise, we conclude that
$$ |\hat g_1(r_1 - r_2)| \leq
\frac{|g'_1(1)-g'_1(0)| + \| g''_1 \|_1}{4\pi^2 |r_1-r_2|^2}
$$
The estimates for $g_2$ and $h_2$ are simpler in nature.
We use the bounds
$$ |\hat g_2(r_1 - r_2)| \leq
\frac{|g_2(1)-g_2(0)| + \| g'_2 \|_1}{2\pi |r_1-r_2|} ,
$$
$$ |\hat h_2(r_1)| \leq
\frac{|h_2(1)-h_2(0)| + \| h'_2 \|_1}{2\pi |r_1|},
$$
as well as
$$|\hat g_3(r_1 - r_2)| \leq \| g_3 \|_1,$$
$$|\hat h_3(r_1)| \leq \| h_3 \|_1,$$
and
$$|\hat \Psi_{\partial_1^2\partial_2^2 F}(r_1,r_2)| \leq
\| \Psi_{\partial_1^2\partial_2^2 F}\|_1.$$
Putting all these together, we see that
$$ |\hat \Psi_F(r_1,r_2)| = O\left(
\frac{1}{|r_1r_2|(r_1-r_2)^2}
+ \frac{1}{|r_1(r_1-r_2)|r_2^2}
+ \frac{1}{r_1^2r_2^2} \right) $$
which is easily verified to be summable over all admissible values of
$r_1$ and $r_2$. We will return to this shortly for a more explicit estimate.
\end{enumerate}

\begin{figure}[t]
\begin{center}
\includegraphics[height=4in]{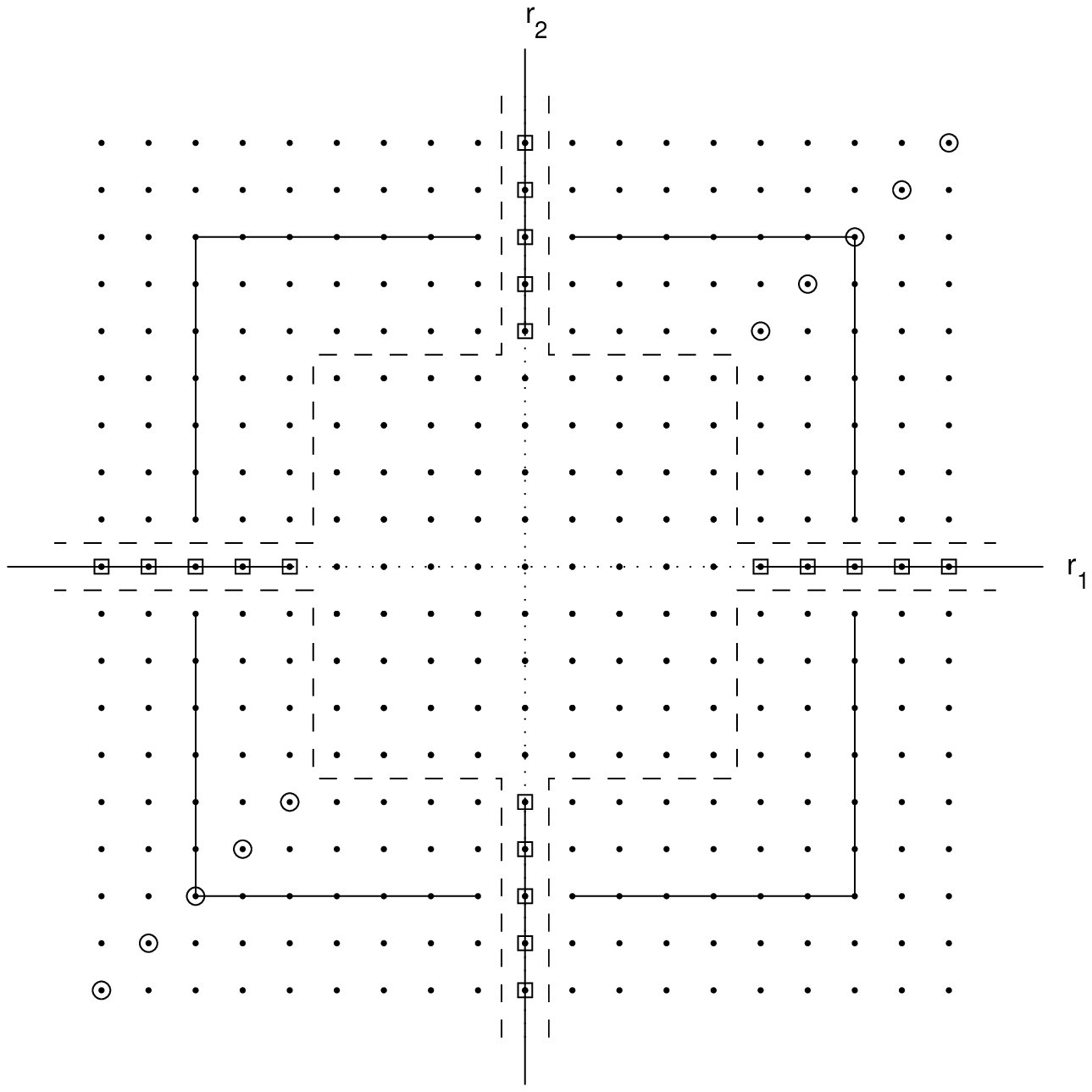}
\end{center}
\caption{\label{figure2} }
\end{figure}

\end{proof}

\subsection*{Explicit numerical estimates}

In this subsection we will work with the specific function $\varphi$
in (\ref{def-varphi}) for which
$$ F(t_1,t_2) = -40(1-t_1)(1-t_2)\left(1-(2-t_1-t_2)^6\right).$$

\subsubsection*{Estimating the value of $C_{\mathrm{axial}}$}
Since $F$ is symmetric, we have
$\hat \varphi(r,0) = \hat \varphi(0,r)$.
We use the formula
(\ref{on-the-axes}) to evaluate $\hat \varphi(r_1,0)$.
It is a simple calculation to show that
$$
J_1(t_1) = -15(1-t_1)+\frac{240}{7}(1-t_1)^2-20(1-t_1)^3+\frac{5}{7}(1-t_1)^9.
$$
Using this expression, we find that $|J_1'(0)-J_1'(1)| = 15$ and
$$\int_0^1 |J_1''| \leq \left(\int_0^1 (J_1'')^2 \right)^{1/2} = 8\sqrt{15}.$$
Hence by (\ref{axial-decay-1}), we obtain the estimate
$$ |\hat \varphi(r,0)| = |\hat \varphi(0, r)| \leq
\frac{15 + 8\sqrt{15}}{4 \pi^2}\frac{1}{r^2}$$
If we define
\beq
C_{\text{axial}}(N) =  \sum_{|r| \leq N}
\left (|\hat \varphi(r,0)| + |\hat \varphi(0, r)| \right ),
\eeq
then it follows that
$$0 \leq  C_{\text{axial}} - C_{\text{axial}}(N)
\leq \frac{15 + 8\sqrt{15}}{\pi^2}\frac{1}{N} <  \frac{5}{N}.$$
To estimate $C_{\text{axial}}(N)$, we still need the actual expression
for $\hat \varphi(r,0)$, which is given in (\ref{r0}).
Taking $N = 50000$, numerical computation shows that
$C_{\text{axial}}(N) = 2.90278\ldots$; hence it follows that
$$
2.90278 \leq C_{\text{axial}} \leq 2.90289.
$$

\subsubsection*{Estimating the value of $C_{\mathrm{main}}$}

We shall estimate the diagonal terms first. We have
$$ g_1(t) = -240t(1-t) $$
from which we obtain
$$ |\hat g_1(0)| = 40, $$
and
\begin{eqnarray}
|\partial_1 \partial_2 F(t_1,t_2)|
& = & \big | 1200(1-t_1)(1-t_2)(2-t_1-t_2)^4 + 280(2-t_1-t_2)^6-40 \big |,
\nonumber \\
& \leq & 1200(1-t_1)(1-t_2)(2-t_1-t_2)^4 + 280(2-t_1-t_2)^6+40
\nonumber
\end{eqnarray}
from which we obtain
$$\| \Psi_{\partial_1 \partial_2 F} \|_{L^1([0,1)^2)} \leq 80.
$$
Hence by (\ref{diagonal-decay-1}), we obtain the estimate
$$ |\hat \varphi(r,r)| \leq \frac{30}{\pi^2}\frac{1}{r^2}$$
It follows that
$$ \sum_{|r| = N+1}^\infty |\hat \varphi(r,r)|
\leq \frac{60}{\pi^2}\frac{1}{N} $$

We next estimate $\hat \varphi(r_1,r_2)$ in the case when
$r_1 \not= r_2$. We begin by noting that
$$
\hat g_1(r_1-r_2) = \frac{120}{\pi^2 (r_1-r_2)^2},\qquad \qquad r_1 \not= r_2.
$$
We have
\begin{eqnarray}
 g_2(t) & = & 240(1+5t(1-t)) \nonumber \\
 h_2(t) & = & -40+280(1-t)^6  \nonumber \\
 g_3(t) & = & 240(-12-15t+20t^2)  \nonumber \\
 h_3(t) & = & -1680(1-t)^5  \nonumber \\
 \partial^2_1\partial^2_2F(t_1,t_2) & = & 14400(2-t_1-t_2)^2
(5+t_1^2+t_2^2-5t_1-5t_2+3t_1t_2)  \nonumber
\end{eqnarray}
from which we obtain
\begin{eqnarray}
\hat g_2(r_1-r_2) & = & -\frac{600}{\pi^2 (r_1-r_2)^2},
\qquad \qquad r_1 \not= r_2. \nonumber \\
|\hat h_2(r_1)| & \leq & \frac{280}{\pi |r_1|},
\qquad \qquad r_1 \not= 0. \nonumber \\
|\hat g_3(r_1-r_2)| & \leq & 3080 \nonumber \\
|\hat h_3(r_1)| & \leq & 280 \nonumber \\
\|\Psi_{\partial^2_1\partial^2_2F} \|_1 & = & 2800.\nonumber
\end{eqnarray}
Putting these together, we obtain the estimate
\begin{equation}
| \hat \varphi(r_1,r_2) |
\leq \frac{105}{\pi^4}\frac{1}{|r_1 r_2|(r_1-r_2)^2}
+  \frac{420}{\pi^4}\frac{1}{r_1^2 r_2^2}
\end{equation}

The following is a simple lemma:

\begin{lemma}
\beq
\sum_{R = N+1}^\infty
\sum_{\max(|r_1|,|r_2|)=R \atop \min(|r_1|,|r_2|)\not=0}
\frac{1}{r_1^2 r_2^2}
< \frac{4\pi^2}{3} \frac{1}{N}
\eeq

\beq
\sum_{R = N+1}^\infty
\sum_{\max(|r_1|,|r_2|)=R \atop {\min(|r_1|,|r_2|)\not=0 \atop r_1 \not= r_2}}
\frac{1}{|r_1r_2|(r_1-r_2)^2}
< 4\left(\frac{\pi^2}{3} +1\right) \frac{1}{N}
\eeq
\end{lemma}

\begin{proof}
The first inequality simply follows from
$$
\sum_{\max(|r_1|,|r_2|)=R \atop \min(|r_1|,|r_2|)\not=0}
\frac{1}{r_1^2 r_2^2}
= \frac{8}{R^2} \sum_{r=1}^{R-1} \frac{1}{r^2}
\leq \frac{4\pi^2}{3}\frac{1}{R^2}.
$$
For the second inequality, we first use the symmetries to write
$$
\sum_{\max(|r_1|,|r_2|)=R \atop {\min(|r_1|,|r_2|)\not=0 \atop r_1 \not= r_2}}
\frac{1}{|r_1r_2|(r_1-r_2)^2}
=
\frac{4}{R} \sum_{r=1}^{R-1} \frac{1}{r(R-r)^2}
+
\frac{4}{R} \sum_{r=1}^{R} \frac{1}{r(R+r)^2}
- \frac{1}{2R^4}.
$$
Using the identity
$$
\frac{1}{r(R-r)^2} = \frac{1}{Rr(R-r)}+\frac{1}{R(R-r)^2},
$$
and Cauchy-Schwarz inequality
we have
$$
\frac{4}{R} \sum_{r=1}^{R-1} \frac{1}{r(R-r)^2} =
\frac{4}{R^2} \left(\sum_{r=1}^{R-1} \frac{1}{r(R-r)}
+ \sum_{r=1}^{R-1} \frac{1}{(R-r)^2} \right)
< \frac{4\pi^2}{3} \frac{1}{R^2}.
$$
For the remaining terms, we use the trivial estimate
$$
\frac{4}{R} \sum_{r=1}^{R} \frac{1}{r(R+r)^2} - \frac{1}{2R^4}
< \frac{4}{R^2}
$$
Hence
$$
\sum_{\max(|r_1|,|r_2|)=R \atop {\min(|r_1|,|r_2|)\not=0 \atop r_1 \not= r_2}}
\frac{1}{|r_1r_2|(r_1-r_2)^2}
< 4\left(\frac{\pi^2}{3} +1\right)\frac{1}{R^2}
$$
and the result follows.
\end{proof}

If we define
\begin{equation}
C_{\text{main}}(N) = \sum_{R = 1}^N
\sum_{\max(|r_1|,|r_2|)=R \atop \min(|r_1|,|r_2|)\not=0}
|\hat \varphi(r_1,r_2)|,
\end{equation}
then we have

\begin{equation}
0 \leq C_{\text{main}} - C_{\text{main}}(N)
< \left (\frac{340}{\pi^2} + \frac{420}{\pi^4} \right) \frac{1}{N}
< \frac{40}{N}.
\end{equation}

For $N=4000$, numerical computation using the formulas
(\ref{rr}) and (\ref{r1r2}) reveals that
$C_{\text{main}}(N) = 4.75145\ldots$; hence with the above error
estimate, we have
\begin{equation}
4.75145\leq C_{\text{main}} \leq 4.76146.
\end{equation}

\subsection*{Explicit expressions for $\hat \varphi(r_1,r_2)$}

The following formulas have been computed using Mathematica, though
it is also possible to compute them easily
using the iterative procedure based on integration by parts which was
outlined in this section earlier.

\begin{multline}\label{r0}
\hat \varphi(r,0) =
\frac{15}{4\pi^2 r^2}\left( 1 - \frac{6}{\pi^2 r^2}
+\frac{45}{\pi^4 r^4}-\frac{135}{\pi^6 r^6}\right) \\
-i\, \frac{60}{7\pi^3 r^3}\left( 1 + \frac{63}{8\pi^2 r^2}
-\frac{315}{8\pi^4 r^4}+\frac{945}{16 \pi^6 r^6}\right).
\end{multline}

\begin{multline}\label{rr}
\hat \varphi(r,r) =
\frac{10}{\pi^2 r^2}\left( 1 - \frac{21}{\pi^2 r^2}
+\frac{315}{2\pi^4 r^4}-\frac{945}{2\pi^6 r^6}\right) \\
+i\, \frac{55}{\pi^3 r^3}\left( 1 - \frac{126}{11\pi^2 r^2}
+\frac{630}{11\pi^4 r^4}-\frac{945}{11 \pi^6 r^6}\right).
\end{multline}

\begin{multline}\label{r1r2}
\hat \varphi(r,s) = \\
-\frac{1575}{4\, {{\pi }^8}\, {r^6}\, {{(r-s)}^2}}
+\frac{525}{4\,{{\pi }^6}\, {r^4}\, {{(r-s)}^2}}
-\frac{35}{2\, {{\pi }^4}\, {r^2}\, {{(r-s)}^2}}
-\frac{1575}{4\, {{\pi }^8}\,{{(r-s)}^2}\, {s^6}} \\
+\frac{225}{2\, {{\pi }^8}\, r\, {{(r-s)}^2}\, {s^5}}
\frac{525}{4\, {{\pi }^6}\, {{(r-s)}^2}\, {s^4}}
+\frac{225}{2\, {{\pi }^8}\, {r^2}\, {{(r-s)}^2}\,{s^4}}
+\frac{225}{2\, {{\pi }^8}\, {r^3}\, {{(r-s)}^2}\, {s^3}} \\
-\frac{75}{2\, {{\pi }^6}\, r\, {{(r-s)}^2}\,{s^3}}
-\frac{35}{2\, {{\pi }^4}\, {{(r-s)}^2}\, {s^2}}
+\frac{225}{2\, {{\pi }^8}\, {r^4}\, {{(r-s)}^2}\, {s^2}}
-\frac{75}{2\, {{\pi }^6}\, {r^2}\, {{(r-s)}^2}\, {s^2}} \\
+\frac{225}{2\, {{\pi }^8}\, {r^5}\, {{(r-s)}^2}\,s}
-\frac{75}{2\, {{\pi }^6}\, {r^3}\, {{(r-s)}^2}\, s}
+\frac{5}{{{\pi }^4}\, r\, {{(r-s)}^2}\, s} \\
i\, \Big(
-\frac{1575}{4\, {{\pi }^9}\, {r^7}\, {{(r-s)}^2}}
+\frac{525}{2\, {{\pi }^7}\, {r^5}\, {{(r-s)}^2}}
-\frac{105}{2\,{{\pi }^5}\, {r^3}\, {{(r-s)}^2}}
-\frac{1575}{4\, {{\pi }^9}\, {{(r-s)}^2}\, {s^7}} \\
+\frac{225}{2\, {{\pi }^9}\,r\, {{(r-s)}^2}\, {s^6}}
+\frac{525}{2\, {{\pi }^7}\, {{(r-s)}^2}\, {s^5}}
+\frac{225}{2\, {{\pi }^9}\, {r^2}\, {{(r-s)}^2}\, {s^5}}
+\frac{225}{2\, {{\pi }^9}\, {r^3}\, {{(r-s)}^2}\,{s^4}}\\
-\frac{75}{{{\pi }^7}\, r\, {{(r-s)}^2}\, {s^4}}
-\frac{105}{2\, {{\pi }^5}\, {{(r-s)}^2}\, {s^3}}
+\frac{225}{2\,{{\pi }^9}\, {r^4}\, {{(r-s)}^2}\, {s^3}}
-\frac{75}{{{\pi }^7}\, {r^2}\, {{(r-s)}^2}\, {s^3}} \\
+\frac{225}{2\, {{\pi }^9}\, {r^5}\, {{(r-s)}^2}\,{s^2}}
-\frac{75}{{{\pi }^7}\, {r^3}\, {{(r-s)}^2}\, {s^2}}
+\frac{15}{{{\pi }^5}\, r\, {{(r-s)}^2}\, {s^2}}
+\frac{225}{2\,{{\pi }^9}\, {r^6}\, {{(r-s)}^2}\, s} \\
-\frac{75}{{{\pi }^7}\, {r^4}\, {{(r-s)}^2}\, s}
+\frac{15}{{{\pi }^5}\, {r^2}\, {{(r-s)}^2}\, s}\Big)
\end{multline}

\providecommand{\bysame}{\leavevmode\hbox
to3em{\hrulefill}\thinspace}
\providecommand{\MR}{\relax\ifhmode\unskip\space\fi MR }
\providecommand{\MRhref}[2]{%
  \href{http://www.ams.org/mathscinet-getitem?mr=#1}{#2}
} \providecommand{\href}[2]{#2}


\begin{thebibliography}{1}

\bibitem{hofm01}
Gerd Hofmeister, \emph{Thin bases of order two}, J. Number Theory
\textbf{86}
  (2001), no.~1, 118--132.

\bibitem{klot69b}
Walter Klotz, \emph{Eine obere {S}chranke f\"ur die {R}eichweite
einer
  {E}xtremalbasis zweiter {O}rdnung}, J. Reine Angew. Math. \textbf{238}
  (1969), 161--168.

\bibitem{mose60}
L.~Moser, \emph{On the representation of $1,2,\ldots,n$ by sums},
Acta Arith.
  \textbf{6} (1960), 11--13.

\bibitem{mose-poun-ridd69}
L.~Moser, J.~R. Pounder, and J.~Riddell, \emph{{On the cardinality
of $h$-Basis
  for $n$}}, J. London Math. Soc. \textbf{44} (1969), 397--407.

\bibitem{mros79}
A.~Mrose, \emph{{Untere Schranken f{\" u}r die Reichweiten von
Extremalbasen
  fester Ordnung}}, Abh. Math. Sem. Univ. Hamburg \textbf{48} (1979), 118--124.

\bibitem{rohr37a}
H.~Rohrbach, \emph{Ein {B}eitrag zur additiven {Z}ahlentheorie},
Math. Zeit.
  \textbf{42} (1937), 1--30.

\end{thebibliography}
\end{document}